\documentclass[a4paper, 11pt, reqno]{amsart}
\usepackage{cite}
\usepackage{url}
\usepackage{amsmath,amssymb,amsthm,latexsym}
\usepackage{graphicx,epsfig}
\usepackage{mathtools} 

\usepackage{mathrsfs} 

\newcommand{\kommentar}[1]{}

\numberwithin{equation}{section} \numberwithin{lemma}{section}

\title[Twin primes from Riemann zeros]{Twin prime correlations from the pair correlation of Riemann zeros}

\author{J. P. Keating and D. J. Smith}

\address{School of Mathematics, University of Bristol\\
Bristol\\ BS8 1TW\\ UK}

\begin{document}

\maketitle  

\begin{abstract}
We establish, via a formal/heuristic Fourier inversion calculation, that the Hardy-Littlewood twin prime conjecture is equivalent to an asymptotic formula for the two-point correlation function of Riemann zeros at a height $E$ on the critical line.  Previously it was known that the Hardy-Littlewood conjecture implies the pair correlation formula, and we show that the reverse implication also holds.  A smooth form of the Hardy-Littlewood conjecture is obtained by inverting the $E \rightarrow \infty$ limit of the two-point correlation function and the precise form of the conjecture is found by including asymptotically lower order terms in the two-point correlation function formula.
\end{abstract}

\section{Introduction}

The Riemann zeta-function is defined by
\begin{equation}
\label{Dirichlet}
\zeta(s)=\sum_{n=1}^{\infty}\frac{1}{n^s} = \prod_p \left(1 - \frac{1}{p^s}\right)^{-1} 
\end{equation}
when ${\rm Re}(s) > 1$ and then by analytic continuation.  It has zeros at the negative even integers, and infinitely many zeros lying off the real line which are termed the {\it non-trivial zeros}.  The non-trivial zeros lie in the strip $0<{\rm Re}(s)<1$, and the Riemann Hypothesis asserts that they all lie on the {\it critical line} ${\rm Re}(s)=1/2$; that is, they all lie at points $s=1/2+iE$, where $E$ is real.  The {\it functional equation} asserts that the zeta function, when multiplied by $\pi^{-s/2}\Gamma(s/2)$, is symmetric with respect to reflection in the critical line \cite{T}.  The zeros of the zeta function are of central importance in Number Theory.  

There has for a long time been an interest in a conjectural connection between the non-trivial zeros of the Riemann zeta-function and the semiclassical theory of quantum chaotic systems \cite{Be1, Ke1, BeKe, BourKe}.  This is motivated by two observations.  First, the explicit formula relating the zeros and the primes closely resembles the semiclassical  trace formula connecting classical periodic orbits and quantum energy levels in chaotic systems.  Second, Montgomery conjectured that the limiting pair correlation of the zeros coincides with that of the eigenvalues of random complex Hermitian matrices drawn from the Gaussian Unitary Ensemble (GUE) of Random Matrix Theory (RMT) \cite{M}, and proved a theorem consistent with this; and it is one of the central conjectures of Quantum Chaos that the energy levels of generic, classically chaotic, non-time-reversal-symmetric systems should, on the scale of the mean level separation, also exhibit GUE statistics in the semiclassical limit.  

Odlyzko computed statistics close to the $10^{20}$th zero of the zeta function, finding that there is a remarkably close agreement between the numerical data and predictions based on Montgomery's conjecture on the scale of the mean zero spacing, but that there are significant deviations on scales large compared to the mean zero spacing \cite{O}.  In \cite{Be2}, Berry wrote down a formula describing Odlyzko's data uniformly.  This augments the random-matrix limit with lower order terms which describe the pair correlations at a height $E$ on the critical line and which vanish when $E\rightarrow\infty$. 

It was shown in \cite{BK} that the lower order terms identified by Berry may be expressed in terms of the lowest zeros of the zeta function; that is, there is a resurgent relationship between the pair correlation of the high-lying zeros and the positions of the low-lying zeros.  Moreover, an additional contribution from the low-lying zeros not captured in Berry's formula was identified.  For a review of these formulae and a comparison with numerical data, see \cite{BeKe}.  The results described above concerning the pair correlation of the Riemann zeros have been extended to other principal $L$-functions \cite{RS, KS} and to all $n$-tuple correlations \cite{BoKe1, BoKe2, RS}.  For an overview, see \cite{Ke2, Ke3}. 

Connections between RMT and the statistical properties of the Riemann zeta-function and other $L$-functions have subsequently been developed in a number of different directions, including moments \cite{KeSn1, KeSn2, CFKRS}, their connections with divisor correlations \cite{CK1, CK2, CK3, CK4, CK5}, and ratios formulae \cite{CFZ}.  The connections have led to several alternative heuristic derivations of the pair correlation formula with lower order terms first derived for the zeta function in \cite{BK}; see, for example, \cite{BoKe3, BoKe4, CS, CK6, CKNon}.         

In his paper \cite{M}, Montgomery mentioned heuristics based on 
the Hardy-Littlewood conjecture concerning the distribution of prime pairs with a given separation (known as {\it twin primes}) \cite{HL} to justify his conjecture that the limiting pair-correlations of the zeta zeros coincide those of GUE eigenvalues. He did not give the details of this calculation, but it 
has subsequently been repeated with variations several times in the literature, for example in \cite{Ke1, GG, BoKe1, BoKe2}.  Goldston and Montgomery \cite{GM} proved rigorously that the limiting pair correlation conjecture is equivalent to an asymptotic formula for the variance of the number of primes in short intervals, and Montgomery and Soundararajan \cite{MS} proved that this variance formula follows from the Hardy-Littlewood twin prime conjecture, under certain assumptions.  The analysis in \cite{GM} has since been extended to all $L$-functions in the Selberg Class \cite{S, BKS}.  The Hardy-Littlewood twin prime conjecture was also the basis of the original heuristic derivation of the pair correlation formula including lower order terms derived in \cite{BK}; that is, the Hardy-Littlewood twin prime conjecture heuristically implies the the pair correlation formula including all lower order terms.

Our purpose in this brief note is to show via a formal heuristic calculation that the reverse implication is also true, namely that the the pair correlation formula including lower order terms obtained in \cite{BK} can be used to derive the Hardy-Littlewood twin prime conjecture.  Thus, heuristically, the two formulae are equivalent.  We demonstrate this via Fourier inversion of the pair correlation formula.  This calculation extends one reported in \cite{ADD}, where an averaged form of the Hardy-Littlewood twin prime conjecture was obtained from the random-matrix limiting expression for the pair correlations, in a manner similar to \cite{GM}.

We denote the non-trivial zeros of the Riemann zeta-function by $1/2+iE_n$ and assume the Riemann Hypothesis, so that $E_n\in{\mathbb R}$ $\forall n$.  Let $d(E)=\sum_n\delta(E-E_n)$ denote the density of zeros at height $E$ on the critical line.  The pair correlation function may then be defined by
\begin{equation}
R_2(\epsilon) = \left<d(E - \epsilon / 2) d(E + \epsilon / 2)\right>,
\end{equation}
where the angular brackets denote an average with respect to $E$.  It may be expressed in the form
\begin{equation}
R_2(\epsilon) = \bar{d}^2(E) + R_2^{(diag)}(\epsilon) + R_2^{(off)}(\epsilon),
\end{equation} 
where
\begin{equation}
\bar{d}(E) = \frac{1}{2 \pi} \ln\left(\frac{E}{2 \pi}\right)
\end{equation}
is the asymptotic mean density of the zeros, $R_2^{(diag)}(\epsilon)$ represents the contributions from the diagonal terms when the pair correlation function is expressed as a sum over pairs of primes using the explicit formula, and  $R_2^{(off)}(\epsilon)$ represents the off-diagonal terms \cite{Ke1}.

The diagonal part of the correlation function may be evaluated to give
\begin{equation}
\lim_{E\to\infty}\frac{1}{\bar{d}^2(E)} R_2^{(diag)} \left(\epsilon/\bar{d}(E)\right) = - \frac{1}{2 \pi^2 \epsilon^2}.
\end{equation}
Montgomery's conjecture is then equivalent to
\begin{equation}
\label{RMTlimit}
\lim_{E\to\infty}\frac{1}{\bar{d}^2(E)} R_2^{(off)} \left(\epsilon/\bar{d}(E)\right) = \frac{\cos (2 \pi \epsilon)}{2 \pi^2 \epsilon^2}.
\end{equation}

For finite $E$ the diagonal part is given by \cite{BK, BeKe}
\begin{equation}
- \frac{1}{4 \pi^2} \left(\frac{d^2}{dw^2} \ln \zeta(1 + i w) + \sum_{\substack{p \\ k \geq 0}} (\ln p)^2 k p^{- (1 + i \epsilon) (k + 1)}\right)_{w = \epsilon} + \textrm{ c.c.}
\end{equation}

The off-diagonal correlations of Riemann zeros may be expressed in terms of primes through the von Mangoldt function;
\begin{equation}
\Lambda(n) = \begin{dcases*} \log p & if $n$ is a power of a prime, $p,$ \\ 0 & otherwise; \end{dcases*}
\end{equation}
explicitly \cite{Ke1},
\begin{equation} \label{Off-diagonal sum}
R_2^{(off)}(\epsilon) = \frac{1}{4 \pi^2} \sum_{n = 1}^{\infty} \sum_{h \neq 0} \frac{\Lambda(n + h) \Lambda(n)}{n} \exp(i (\epsilon \ln n - E h / n)) + \textrm{c.c.}
\end{equation}
The finite-$E$ expression for the off-diagonal contribution to the pair correlation of Riemann zeros can then be calculated by making use of a conjecture by Hardy and Littlewood \cite{HL} concerning the distribution of pairs of primes. This states that if 
\begin{equation}
\alpha(h) = \lim_{N \to \infty} \frac{1}{N} \sum_{n = 1}^N \Lambda(n + h) \Lambda(n),
\end{equation}
then
\begin{equation}
\alpha(h) = \begin{dcases*} 2 C_2 \prod_{p | h} \frac{p - 1}{p - 2} & if $h$ is even \\ 0 & otherwise, \end{dcases*}
\end{equation}
where $C_2$ is the twin prime constant
\begin{equation}
C_2 = \prod_{p > 2} \left(1 - \frac{1}{(p - 1)^2}\right) \approx 0.6601618.
\end{equation}

Using the Hardy-Littlewood conjecture, the finite-$E$ off-diagonal contribution to the two-point correlation function was evaluated in \cite{BK} to give
\begin{equation}
\label{largeE}
R_2^{(off)}(\epsilon) = \frac{1}{4 \pi^2} |\zeta(1 + i \epsilon)|^2 \left(\exp(- 2 \pi i \epsilon \bar{d}(E)) \prod_p \left(1 - \left(\frac{1 - p^{- i \epsilon}}{p - 1}\right)^2\right)\right) + \textrm{ c.c.}
\end{equation}
This formula has since been derived using a number of other approaches \cite{BoKe3, BoKe4, CS, CK6, CKNon}.  It reduces to \eqref{RMTlimit} in the limit $E\to\infty$.

The purpose of this note is to demonstrate a dual relationship between the correlations of zeros and primes. A smoothed form of the Hardy-Littlewood conjecture will first be determined by making use of the limiting form of the pair correlation of zeros \eqref{RMTlimit} and then the full Hardy-Littlewood conjecture will be recovered using finite $E$ corrections to the pair correlation, as captured by \eqref{largeE}.

\section{Inverting the 2-point correlation}
The off-diagonal contribution \eqref{Off-diagonal sum} can be expressed as
\begin{equation}
R_2^{(off)}(\epsilon) = \frac{1}{4 \pi^2} \int_1^{\infty} \int_{- \infty}^{\infty} \frac{\alpha(h)}{y} \exp(i (\epsilon \ln y - E h / y)) \, dh dy + \textrm{ c.c,}
\end{equation}
c.f.~\cite{Ke1, BoKe1, BoKe2}, and so an expression for $\alpha(h)$ may be determined by Fourier inversion, as follows.

\begin{align}
\int_{- \infty}^{\infty} R_2^{(off)}(\epsilon) \, d\epsilon & = \frac{1}{2 \pi} \int_1^{\infty} \int_{- \infty}^{\infty} \frac{\alpha(h)}{y} \delta(\ln y) \exp(-  i E h / y) \, dh dy + \textrm{ c.c }\\ 
& = \frac{1}{2 \pi} \int_0^{\infty} \int_{- \infty}^{\infty} \alpha(h) \delta(z) \exp(- i E h \exp(- z)) \, dh dz + \textrm{ c.c }\\
& = \frac{1}{2 \pi} \int_{- \infty}^{\infty} \alpha(h) \exp(- i E h) \, dh.
\end{align}
Giving the Fourier transform 
\begin{equation}
\mathscr{F}[\alpha(h)] = 2 \pi \int_{- \infty}^{\infty} R_2^{(off)}(\epsilon) \, d\epsilon,
\end{equation}
and by Fourier inversion,
\begin{equation}
\label{inv}
\alpha(h) = \int_{- \infty}^{\infty} \int_{- \infty}^{\infty} R_2^{(off)}(\epsilon) \exp(i h E) \, d\epsilon dE.
\end{equation}

\section{Averaged form of the Hardy-Littlewood \\ conjecture}
In \cite{Ke1} Montgomery's expression for the pair correlation was shown heuristically to follow from knowledge of the average of the function in the Hardy-Littlewood conjecture. Asymptotically as $|h| \to \infty$ the average is given by $1 - \frac{1}{2 |h|},$ in the sense that
\begin{equation}
\frac{1}{2 |h|} \sum_{H = - h}^h \alpha(H) \sim 1 - \frac{\ln |h|}{2 |h|} \textrm{ as } |h| \to \infty.
\end{equation}
The $- \frac{1}{2 |h|}$ term in the averaged Hardy-Littlewood conjecture was recovered in \cite{ADD} by assuming Montgomery's conjecture (the first term does not contribute to the zero correlations), specifically by inverting the off-diagonal contribution of the limiting pair correlation \eqref{RMTlimit}. For completeness, this calculation is shown here to follow from the inversion formula \eqref{inv}.

Substituting the limiting form of the off-diagonal pair correlation gives the integral over $\epsilon,$
\begin{align}
\frac{1}{2 \pi^2} \int_{- \infty}^{\infty} \frac{1}{\epsilon^2} \cos(2 \pi \epsilon \bar{d}(E)) \, d\epsilon = & \frac{1}{\pi} \int_{- \infty}^{\infty} \frac{1}{\epsilon^2} \cos(\epsilon \bar{d}(E)) \, d\epsilon \nonumber \\
= & \frac{1}{\pi} \mathscr{F}\left[\frac{1}{\epsilon^2}\right](\bar{d}(E)).
\end{align}

This Fourier transform may be deduced by considering the Fourier transform of the triangle function
\begin{equation}
T(x) = \begin{dcases*} 1 - |x| & if $|x| \leq 1$ \\ 0 & otherwise, \end{dcases*}
\end{equation}
and by noting that the triangle function may be expressed in terms of the sign function,
\begin{equation}
\textrm{sgn}(x) = \begin{dcases*} - 1 & if $x \leq 0$ \\ 1 & otherwise, \end{dcases*}
\end{equation}
by the linear relation
\begin{equation}
T(x) = \frac{1}{2} (1 - x) \textrm{sgn}(1 - x) + \frac{1}{2} (1 + x) \textrm{sgn}(1 + x) - x \textrm{sgn}(x).
\end{equation}
Then 
\begin{equation}
\mathscr{F}[T](k) = 2 \int_0^1 (1 - x) \cos(k x) \, dx = \left(\textrm{sinc}\left(\frac{k}{2}\right)\right)^2,
\end{equation}
by integrating by parts, giving
\begin{align}
\int_{- \infty}^{\infty} \left(\textrm{sinc}\left(\frac{k}{2}\right)\right)^2 \textrm{exp}(i k x) \, dk = & \pi (1 - x) \textrm{sgn}(1 - x) \nonumber \\
& +  \pi (1 + x) \textrm{sgn}(1 + x) \nonumber \\
& - 2 \pi x \textrm{sgn}(x)
\end{align}
by Fourier inversion. Expanding the sinc function in terms of exponential functions then gives
\begin{align}
\int_{- \infty}^{\infty} \left(\textrm{sinc}\left(\frac{k}{2}\right)\right)^2 \textrm{exp}(i k x) \, dk = & - \mathscr{F}\left[\frac{1}{k^2}\right](1 - x) \nonumber \\
& - \mathscr{F}\left[\frac{1}{k^2}\right](1 + x) \nonumber \\
& + 2 \mathscr{F}\left[\frac{1}{k^2}\right](x),
\end{align}
yielding the identification
\begin{equation}
\mathscr{F}\left[\frac{1}{x^2}\right](k) = - \pi k \textrm{sgn}(k).
\end{equation}
The smoothed correlation function is then given by
\begin{align}
& - \mathscr{F} \left[\bar{d}(E) \textrm{sgn}(\bar{d}(E))\right](h) \nonumber \\
& \sim \frac{1}{\pi} \int_0^1 \ln (E) \cos (h E) \, dE \textrm{ as } |h| \to \infty \nonumber \\
& \sim - \frac{1}{\pi h} \textrm{Si}(h) \textrm{ as } |h| \to \infty.
\end{align}
This then gives $- \frac{1}{2 |h|}$ as $\int_0^{\infty} \frac{\sin x}{x} \, dx = \frac{\pi}{2};$ recovering the result found in \cite{ADD}. 

\section{The full Hardy-Littlewood conjecture}
\subsection{Off-diagonal pair correlation}
As reviewed in the introduction, various heuristic calculations \cite{BK, BoKe3, CS, CK6, CKNon} lead to an expression for the off-diagonal contribution to the two-point correlation function;
\begin{equation}
\frac{1}{2 \pi^2} \exp(- 2 \pi i \epsilon \bar{d}(E)) |\zeta(1 + i \epsilon)|^2 \prod_p \frac{\left(1 - \frac{1}{p^{1 + i \epsilon}}\right) (p^2 - 2 p + p^{1 - i \epsilon})}{(p - 1)^2}.
\end{equation}
Using the Euler product on the $1$-line allows us to rewrite this as
\begin{align}
& & \frac{1}{2 \pi^2} \exp(- 2 \pi i \epsilon \bar{d}(E)) \zeta(1 - i \epsilon) \prod_p \left(1 - \frac{1}{(p - 1)^2} + \frac{p^{1 - i \epsilon}}{(p - 1)^2}\right) \nonumber \\
& = & \frac{1}{2 \pi^2} \exp(- 2 \pi i \epsilon \bar{d}(E)) \zeta(1 - i \epsilon) \prod_p \left(1 + \frac{p^{1 - i \epsilon} \left(1 - \frac{1}{p^{1 - i \epsilon}}\right)}{(p - 1)^2}\right).
\end{align}

Euler's totient function $\phi(n)$ is the number of natural numbers less than and coprime to $n$; and the M$\ddot{\textrm{o}}$bius function $\mu(n)$ is zero if $n$ is divisible by a square, and $-1$ to the power  of the number of distinct prime factors of n otherwise \cite{HW}.  Both are multiplicative functions.  This allows the products over primes to be changed to sums over natural numbers in the following way.
\begin{equation}
\frac{1}{2 \pi^2} \exp(- 2 \pi i \epsilon \bar{d}(E)) \zeta(1 - i \epsilon) \prod_p \left(1 + \left(\frac{\mu(p)}{\phi(p)}\right)^2 p^{1 - i \epsilon} \left(1 + \frac{\mu(p)}{p^{1 - i \epsilon}}\right)\right),
\end{equation}
which, by multiplicativity of $\phi$ and $\mu$ yields
\begin{equation}
\frac{1}{2 \pi^2} \exp(- 2 \pi i \epsilon \bar{d}(E)) \zeta(1 - i \epsilon) \sum_{n = 1}^{\infty} \left(\frac{\mu(n)}{\phi(n)}\right)^2 n^{1 - i \epsilon} \prod_{p | n} \left(1 + \frac{\mu(p)}{p^{1 - i \epsilon}}\right),
\end{equation}
where the sum is over squarefree values of $n$ due to the presence of $\mu(n).$ Making use of the multiplicativity of $\mu$ once more gives
\begin{equation}
\frac{1}{2 \pi^2} \exp(- 2 \pi i \epsilon \bar{d}(E)) \zeta(1 - i \epsilon) \sum_{n = 1}^{\infty} \left(\frac{\mu(n)}{\phi(n)}\right)^2 n^{1 - i \epsilon} \sum_{d | n} \frac{\mu(d)}{d^{1 - i \epsilon}}.
\end{equation}
Using now the series representation of the Riemann zeta-function on the $1$-line results in
\begin{equation}
\frac{1}{2 \pi^2} \exp(- 2 \pi i \epsilon \bar{d}(E)) \sum_{n = 1}^{\infty} \left(\frac{\mu(n)}{\phi(n)}\right)^2 n^{1 - i \epsilon} \sum_{m = 1}^{\infty} \sum_{d | n} \frac{\mu(d)}{(md)^{1 - i \epsilon}}.
\end{equation}
Replacing $m d$ by $l,$ the sum over $m$ and $d | n$ becomes a sum over $l,$ $d | n$ and $d | l,$ 
\begin{equation}
\sum_{l = 1}^{\infty} \sum_{d | n, d | l} \frac{\mu(d)}{l^{1 - i \epsilon}} = \sum_{l = 1}^{\infty}  \frac{1}{l^{1 - i \epsilon}} \sum_{d | (l, n)} \mu(d),
\end{equation}
which, by the indicator property of the the M$\ddot{\textrm{o}}$bius function \cite{HW} gives
\begin{equation}
 \sum_{\substack{l =1 \\ (l,n) = 1}}^{\infty} \frac{1}{l^{1 - i \epsilon}}.
\end{equation}
The two-point correlation function is therefore given by
\begin{equation}
\frac{1}{2 \pi^2} \exp(- 2 \pi i \epsilon \bar{d}(E)) \sum_{n = 1}^{\infty} \sum_{\substack{l =1 \\ (l,n) = 1}}^{\infty} \left(\frac{\mu(n)}{\phi(n)}\right)^2 \left(\frac{n}{l}\right)^{1 - i \epsilon},
\end{equation}
which may be approximated by
\begin{equation}
\frac{1}{2 \pi^2} \exp(- 2 \pi i \epsilon \bar{d}(E)) \sum_{n = 1}^{\infty} \sum_{\substack{l =1 \\ (l,n) = 1}}^{\infty} \left(\frac{\mu(n)}{\phi(n)}\right)^2 \int_0^{E / 2 \pi} \frac{1}{y^{1 - i \epsilon}} \delta(y - l / n) \, dy,
\end{equation}
where the resulting error is reduced by taking the value of $E$ to be large. Substituting the expression for the mean density and taking the exponential inside the integral then gives
\begin{equation}
\frac{1}{2 \pi^2} \sum_{n = 1}^{\infty} \sum_{\substack{l =1 \\ (l,n) = 1}}^{\infty} \left(\frac{\mu(n)}{\phi(n)}\right)^2 \int_0^1 \frac{1}{z} \exp(i \epsilon \ln z) \delta(E / (2 \pi) z- l / n) \, dz,
\end{equation}
after making a change of variables in the integral.

\subsection{Full inversion}
Substituting this into the Fourier inversion formula gives 

\begin{align}
& \frac{1}{2 \pi^2} \sum_{n = 1}^{\infty} \sum_{\substack{l =1 \\ (l,n) = 1}}^{\infty} \left(\frac{\mu(n)}{\phi(n)}\right)^2 \int_{- \infty}^{\infty} \int_{- \infty}^{\infty} \int_0^1 \frac{1}{z} \exp\left(i \epsilon \ln z\right) \delta(E / (2 \pi) z - l / n) \, dz \nonumber \\
& \exp(i h E) \, d\epsilon dE.
\end{align}
Integrating over $\epsilon$ produces a delta function,
\begin{align}
& \frac{1}{\pi} \sum_{n = 1}^{\infty} \sum_{\substack{l =1 \\ (l,n) = 1}}^{\infty} \left(\frac{\mu(n)}{\phi(n)}\right)^2 \int_{- \infty}^{\infty} \int_0^1 \frac{1}{z} \delta(E / (2 \pi) z - l / n) \delta(\ln z) \nonumber \\
& \exp(i h E) \, dz dE \nonumber \\
= & \frac{1}{\pi} \sum_{n = 1}^{\infty} \sum_{\substack{l =1 \\ (l,n) = 1}}^{\infty} \left(\frac{\mu(n)}{\phi(n)}\right)^2 \int_{- \infty}^{\infty} \int_{- \infty}^0 \delta(w) \delta(E / (2 \pi) \exp(w) - l / n) \nonumber \\
& \exp(i h E) \, dw dE.
\end{align}
Integrating over $w \leq 0$ captures half the mass of the delta function, giving
\begin{align}
& \frac{1}{2 \pi} \sum_{n = 1}^{\infty} \sum_{\substack{l =1 \\ (l,n) = 1}}^{\infty} \left(\frac{\mu(n)}{\phi(n)}\right)^2 \int_{- \infty}^{\infty} \delta(E / (2 \pi) - l / n) \exp(i h E) \, dE \nonumber \\
= & \sum_{n = 1}^{\infty} \sum_{\substack{l =1 \\ (l,n) = 1}}^{\infty} \left(\frac{\mu(n)}{\phi(n)}\right)^2 \int_{- \infty}^{\infty} \delta(u - l / n) \exp(2 \pi i h u) \, du \nonumber \\
= & \sum_{n = 1}^{\infty} \sum_{\substack{l =1 \\ (l,n) = 1}}^{\infty} \left(\frac{\mu(n)}{\phi(n)}\right)^2 \exp((2 \pi i l h) / n),
\end{align}
which may be approximated by
\begin{equation}
\sum_{n = 1}^{\infty} \left(\frac{\mu(n)}{\phi(n)}\right)^2 c_n(h),
\end{equation}
where
\begin{equation}
c_n(h) = \sum_{\substack{l =1 \\ (l,n) = 1}}^n \exp((2 \pi i l h) / n)
\end{equation}
is Ramanujan's sum \cite{HW}.

Ramanujan sums are multiplicative,
\begin{equation}
c_m(l) c_n(l) = c_{mn}(l)
\end{equation}
if $(m, n) = 1.$ They also satisfy
\begin{equation}
c_p(n) = \begin{dcases*} p - 1 & if $p | n$ \\ - 1 & otherwise \end{dcases*}
\end{equation}
for all primes, $p.$ Using these properties then gives
\begin{align}
& \sum_{n = 1}^{\infty} \left(\frac{\mu(n)}{\phi(n)}\right)^2 c_n(h) \nonumber \\
= & \prod_p \left(1 + \left(\frac{\mu(p)}{\phi(p)}\right)^2 c_p(h)\right) \nonumber \\
= & \begin{dcases*} 2 \prod_{p > 2} \left(1 + \left(\frac{\mu(p)}{\phi(p)}\right)^2 c_p(h)\right) & if $h$ is even \\ 0 & otherwise, \end{dcases*}
\end{align}
where 
\begin{align}
& \prod_{p > 2} \left(1 + \left(\frac{\mu(p)}{\phi(p)}\right)^2 c_p(h)\right) \nonumber \\
= & C_2 \prod_{p | h} \left(\frac{p - 1}{p - 2}\right).
\end{align}

Thus the exact Hardy-Littlewood twin prime conjecture may be derived formally/heuristically from the pair correlation formula for the Riemann zeros including lower order terms.  This is what we set out to demonstrate.  As already shown in \cite{BK} the reverse implication also holds heuristically.  Therefore the two conjectures may be thought of as being formally equivalent.

\section{Acknowledgements}
We are grateful to Brian Conrey and Nina Snaith for a conversation which prompted us to perform this calculation.  JPK is pleased to acknowledge support from a Royal Society Wolfson Research Merit Award and ERC Advanced Grant 740900 (LogCorRM).

\end{document}